\documentclass[12pt,reqno]{amsart}
\usepackage[foot]{amsaddr}
\usepackage[utf8]{inputenc}
\usepackage[english]{babel}
\usepackage{amsmath,amssymb,amsfonts,amsbsy,amsthm,amscd,latexsym,graphicx}
\usepackage{empheq,textcomp}
\usepackage{framed}

\theoremstyle{definition}
\newtheorem{theorem}{Theorem} [section]
\newtheorem{corollary}[theorem]{Corollary}
\newtheorem{lemma}[theorem]{Lemma}
\newtheorem{proposition}[theorem]{Proposition}
\newtheorem{definition}[theorem]{Definition}

\newtheorem{remark}[theorem]{Remark}

\newtheorem{example}[theorem]{Example}

\numberwithin{equation}{section}

\newcommand{\C}{\mathbb{C}}

\newcommand{\Fc}{{\mathcal{F}}}

\newcommand{\N}{\mathbb{N}}
\newcommand{\R}{\mathbb{R}}

\newcommand{\Z}{\mathbb{Z}}

\newcommand{\Eq}{\, = \,}

\newcommand{\Le}{\, \le \,}
\newcommand{\Les}{\, \lesssim \,}
\newcommand{\Ge}{\, \geq \,}

\newcommand{\qeddef}{{\quad $\diamondsuit$}}

\newcommand{\bigabs}[1]{\bigl|\,#1\,\bigr|}
\newcommand{\Bigabs}[1]{\Bigl|\,#1\,\Bigr|}

\newcommand{\ip}[2]{\langle\,#1,#2\,\rangle}
\newcommand{\bigip}[2]{\bigl\langle \,#1, \, #2 \,\bigr\rangle}

\newcommand{\norm}[1]{\|\,#1\,\|}
\newcommand{\bignorm}[1]{\bigl\|\,#1\,\bigr\|}
\newcommand{\Bignorm}[1]{\Bigl\|\,#1\,\Bigr\|}
\newcommand{\biggnorm}[1]{\biggl\|\,#1\,\biggr\|}
\newcommand{\bigparen}[1]{\bigl(\,#1\,\bigr)}
\newcommand{\Bigparen}[1]{\Bigl(\,#1\,\Bigr)}

\newcommand{\set}[1]{\{#1\}}
\newcommand{\bigset}[1]{\bigl\{#1\bigr\}}

\newcommand{\clspan}{{\overline{\text{span}}}}

\newcommand{\inN}{_{n\in\N}}
\newcommand{\sumli}{\sum_{n=1}^\infty}

\newcommand{\schtr}{\set{T_{\lambda_n}g}_{\inN}}

\newcommand{\Chi}{\chi_{_{[0,\alpha]}}}
\title{Existence of Unconditional frames of translates in modulation spaces}

\author{Pu-Ting Yu}

\address{Department of Mathematics, University of Oregon, Eugene, OR, 97403, USA}

\email{putingyu@uoregon.edu}
\date{February 2025}

\begin{document}
\begin{abstract}
Let $1\leq p\leq 2$ and let $\Lambda = \set{\lambda_n}\inN \subseteq \R$ be an arbitrary subset. We prove that for any $g\in M^p(\R)$ with $1\leq p\leq 2$ the system of translates $\set{g(x-\lambda_n)}\inN$ is never an unconditional basis for $M^q(\R)$ for $p\leq q\leq p'$, where $p'$ is the conjugate exponent of $p.$ In particular, $M^1(\R)$ does not admit any Schauder basis formed by a system of translates.
We also prove that for any $g\in M^p(\R)$ with $1< p\leq 2$ the system of translates $\set{g(x-\lambda_n)}\inN$ is never an unconditional frame for $M^p(\R).$ Several partial results regarding the existence of unconditional frames formed by a system of translates in $M^1(\R)$ as well as in $M^p(\R)$ with $2<p<\infty$ will be presented as well. 
\end{abstract}

\maketitle
\section{Introduction}
Let $X$ be a separable Banach space. A sequence $\set{x_n}\inN \subseteq X$ is a \emph{Schauder basis} for $X$ if 
there exists a sequence $\set{x_n^*}\inN\subseteq X^*$ satisfying $\ip{x}{x_n^*}=\delta_{nm}$ such that every $x\in X$ can be uniquely expressed as \begin{equation}
\label{basis_expansion}
    x=\sumli \ip{x}{x_n^*}x_n,
\end{equation} with the convergence of the series in the norm of $X.$ Here we use $\ip{x}{x_n^*}$ to mean $x_n^*(x)$. The sequence $\set{x_n^*}\inN$ is called the \emph{coefficient functionals} associated with $\set{x_n}\inN.$ 
We say that a Schauder basis $\set{x_n}\inN$ is an \emph{unconditional basis} for X if the series in Equation $(\ref{basis_expansion})$ converges unconditionally for all $x\in X$, i.e., $\sumli c_{\sigma(n)}x_{\sigma(n)}$ converges for any permutation $\sigma:\N\rightarrow \N$. A major drawback of Schauder bases is exactly the uniqueness of the expression in Equation (\ref{basis_expansion}), which makes an explicit construction of a Schauder basis possessing certain properties arduous in general. A generalized notion of Schauder bases is the \emph{Schauder frames}. Roughly speaking, the idea of Schauder frames is to give up the uniqueness of the representation in Equation (\ref{basis_expansion}) but in return gain more flexibility in the choice of coefficients of the expansions in Equation (\ref{basis_expansion}). We say that a sequence $\set{x_n}\inN\subseteq X$ is a \emph{Schauder frame} if there exists a sequence $\set{x_n^*}\inN\subseteq X^*,$ which is also called the \emph{associated coefficient functionals}, such that Equation (\ref{basis_expansion}) holds with the convergence of the series in norm of $X$ for all $x\in X.$ Similarly, a Schauder frame is called an \emph{unconditional frame} if the series in Equation (\ref{basis_expansion}) converges unconditionally in the norm of $X$ for all $x\in X.$ In general, a Schauder frame that is not a Schauder basis possesses more than one sequence of coefficient functionals (see \cite[Corollary 5.14]{PY24}). Consequently, we do get more flexibility in choice of expressions for an element in $X$. For more background knowledge along this direction as well as other generalizations of Schauder bases, we refer to \cite{AK06} \cite{CHL99}, and \cite{Gro91}.

 The existence of a Schauder basis or a Schauder frame of a specific form for certain Banach spaces has been one of the major research questions in the field of Banach space theory, harmonic analysis approximation theory, etc, in the past decades. For example, Olson and Zalik conjectured in \cite{OZ92} that for any $g\in L^2(\R)$ and any subset $\set{\lambda_n}\inN\subseteq \R$ the sequence of translates $\set{g(x-\lambda_n)}\inN$ can never be a Schauder basis for $L^2(\R).$ This conjecture remains open as of the time of writing. In the same paper where Olson-Zalik conjecture was formulated, Olson and Zalik proved that such a sequence can never be an unconditional basis for $L^2(\R)$. By proving certain density results associated with irregular Gabor systems, Christensen, Deng and Heil showed in \cite{CDH99} that $L^2(\R)$ does not admit a sequence of translates satisfying the inequality $A\norm{f}^2_{L^2(\R)}\leq \sumli|\ip{f}{T_{\lambda_n}g}|^2\leq B\norm{f}^2_{L^2(\R)}$  for all $f\in L^2(\R)$ whatever positive constants $A$ and $B$ are. The existence of unconditional frames formed by systems of translates in $L^2(\R)$ was disproved by Lev and Tselishcehv in \cite{LT25} (In fact, they proved no sequence of translates can be an unconditional frame for $L^p(\R)$ for any $1\leq p\leq 2$). Currently, the best partial result known is due to Odell et al., who showed in \cite{OSSZ11} that if $g\in L^1(\R)\cap L^p(\R)$, then $\set{g(x-\lambda_n)}\inN$ can never be a Schauder frame for $L^p(\R)$ for any $1<p<\infty$ if $\set{\lambda_n}\inN$ is uniformly discrete (see Lemma \ref{Schauder_discre_index} for the definition). It is worth noting that a stronger conjecture regarding a certain density result of Gabor Schauder basis was formulated by Deng and Heil in \cite{DH00}. In particular, a positive answer to the Deng-Heil conjecture yields a positive answer to the Olson-Zalik conjecture. Aside from $L^2(\R)$, the Olson-Zalik conjecture has been analogously formulated and studied in the setting of $L^p(\R)$ for any $1\leq p<\infty$. For relevant research papers, we refer to \cite{FOSZ14}, \cite{LT25}, \cite{LT25b}, \cite{LT25c} and \cite{OSSZ11} (see also \cite{GMZ16} for more open problems related to Schauder bases in Banach spaces).

 While it is known that $L^2(\R)$ does not admit any unconditional basis or unconditional frames formed by a sequence of translates, the complete characterization of closed subspaces of $L^2(\R)$ that do not admit an unconditional basis (or an unconditional frame) formed by a sequence of translates is still unknown. We accordingly propose the following question: 
 ``\emph{Assume that $X$ is a Banach space that is continuously embedded into $L^2(\R)$ and is invariant under Fourier transform. Is it true that $X$ does not admit an unconditional basis or even does not admit an unconditional frame formed by a sequence of translates}?" 
 Among those Banach spaces that are invariant under Fourier transform and are continuously embedded into $L^2(\R)$, there is a family of such Banach spaces of our particular interest due to its pivotal role in the field of time-frequency analysis -- \emph{modulation spaces}. Since it was discovered in the early 1980s by Feichtinger and then developed in a series of joint work with Gr\"{o}chenig, modulation spaces have been recognized as the proper spaces for time-frequency analysis. For textbook recountings on time-frequency analysis, we refer to \cite{BO20} and \cite{Gro01}. It is known that modulation spaces, denoted by $M^p(\R)$, are invariant under Fourier transform and are continuously embedded into $L^2(\R)$ for any $1\leq p \leq 2$. In this paper, we study the existence of unconditional bases and unconditional frames formed by a system of translates in $M^p(\R).$ Specifically, our results can be summarized as follows. 
\begin{enumerate}
\setlength\itemsep{0.2em}
    \item [\textup{(a)}] There does not exist an unconditional basis formed by a sequence of translates for $M^p(\R)$ for any $1\leq p\leq 2.$ In particular, $M^1(\R)$ does not admit any Schaduer basis formed by a system of translates. For the case $p>2$ we prove that 
     the system of translates $\set{g(x-\lambda_n)}\inN$ is never an unconditional basis if $g\in L^{p'}(\R)\cap M^p(\R)$ whatever $\set{\lambda_n}\inN\subseteq \R $ is.
    \item [\textup{(b)}] There does not exist an unconditional frame formed by a sequence of translates for $M^p(\R)$ for any $1< p\leq 2.$ 
    However, there do exist unconditional frames formed by a sequence of translates for $M^p(\R)$ for any $2<p<\infty.$
\end{enumerate} 
 

 This paper is organized as follows. Section 2 will be devoted to notations, definitions and some known results required for this paper. We then present our main results in Section 3. Specially, we will start with the case $1<p\leq 2$ and then close this paper by the case $p=1.$ Partial results regarding the case $p>2$ will be presented as well.
 \section{Preliminaries}
\label{preliminaries}
Throughout this paper, an element with a star super index, such as $x^*$, means an element in the dual space $X^*$ of $X$ whenever $X$ is the Banach space that is specifically described in the corresponding statement.  For the sake of notational consistence, we will use $\ip{x}{x^*}$ to denote $x^*(x)$ whenever $x\in X$ and $x^*\in X^*.$ We will use $p'$ to denote the conjugate exponent of $p$, i.e., $\frac{1}{p}+\frac{1}{p'}=1$ when $1\leq p\leq \infty.$ Finally, a sequence/system of translates in a Banach space $X$ is a sequence of the form $\set{T_{\lambda_n}g}\inN$ for some $g\in X$ and some subset $\set{\lambda_n}\inN$ of $\R.$ Here $T_\lambda$ denotes the translation operator defined (on an appropriate domain) by $(T_\lambda g)(x)=g(x-\lambda).$

We define \emph{modulation spaces} as follows.
\begin{definition}\label{Modulation_definition}
Fix a Schwartz function $\psi\in S(\R)$. 
\begin{enumerate}\setlength\itemsep{0.5em}
    \item [\textup{(a)}] Let $f\in S'(\R)$ be a tempered distribution. The \emph{short-time Fourier transform} of $f$, denoted by $V_{\psi}f$, is the measurable function on $\R^2$ defined by  $$V_{\psi}f(x,w)\Eq \overline{\ip{M_wT_x \psi}{f}} \Eq  \ip{f}{M_wT_x \psi}.$$

     \item [\textup{(b)}] For $1\leq p \leq \infty$, the \emph{modulation space} $M^{p}(\R)$ is the space consisting of 
     all tempered distributions $f$ for which $\norm{V_\psi f}_{L^{p}(\R^2)}$ is finite, i.e.,
     $$M^{p}(\R)=\bigset{f \in S'(\R)\,\big|\,  \norm{f}_{M^{p}(\R)}\Eq \norm{V_\psi f}_{L^{p}(\R^2)}<\infty},$$
      with the usual modifications if $p=\infty$ or $q=\infty.$ \qeddef
    \end{enumerate}
\end{definition}

The definition that we adopt in this paper of modulation spaces is in fact a specific family of modulation spaces in literature. Using the mixed 
 $L^{p,q}$-norm $\norm{\cdot}_{L^{p,q}(\R^2)}$, one can generalize $M^p(\R)$ to $M^{p,q}(\R)$ for any $1\leq q\leq \infty$. Following the terminology in \cite{Fei06}, we briefly mention that another further generalization of modulation spaces is the \emph{Wiener amalgam spaces} $W(B,C)$. The modulation spaces we defined above are indeed $W(\Fc(L^p),\ell^p)$ in disguise.
 For more history regarding modulation spaces as well as related generalizations, we refer to \cite{Fei06}.

It is known that $M^p(\R)$ are Banach spaces that are invariant under the Fourier transform. Surprisingly, the definition of $M^p(\R)$ is independent of the choice of the Schwartz function $\psi$ that is used in the definition. That is, using a different $\psi\in S(\R)$ in the definition of modulation spaces will yield an equivalent norm. Next, we collect several results related to modulation spaces in the following theorem. We remark that some of them are combinations of two different results. For notational convenience, we will use $A\lesssim_{\phi}B$ to mean that there exists some constant $C>0$ depending on $\phi$ such that $A\leq CB.$ By the notation $X\hookrightarrow Y$ we mean that $X,Y$ are two Banach spaces satisfying $X\subseteq Y$ and $\norm{\cdot}_{Y}\Le C\norm{\cdot}_{X}$ for some constant $C>0$.
\begin{theorem}
\label{tf_lemmas}
Let $1\leq p,q <\infty.$
\begin{enumerate} \setlength\itemsep{0.5em}

 \item [\textup{(a)}] (\cite[Theorem 11.3.6]{Gro01} and \cite[Theorem 12.2.2]{Gro01}) For any $p\leq q$ we have  $$M^{1}(\R)\stackrel{\text{dense}}{\subseteq} M^{p}(\R)\stackrel{\text{dense}}{\subseteq} M^{q}(\R)\subseteq M^{\infty}(\R)$$
    with $\norm{\cdot}_{M^{p}(\R)}\lesssim_{p,q} \norm{\cdot}_{M^{q}(\R)}$. 
Moreover, we have that $(M^p(\R))^*=M^{p'}(\R).$

 \item [\textup{(b)}] (\cite[Theorem 11.3.5]{Gro01})  Modulation spaces are invariant under Fourier transform. Moreover, we have  $\norm{f}_{M^p(\R)}\lesssim_p \norm{\widehat{f}}_{M^{p}(\R)}$ for all $f\in M^p(\R).$ 

\item [\textup{(c)}](\cite[Theorem 11.3.6]{Gro01} and \cite[Remark 3.2]{CLS11})  Assume that $\set{T_{\lambda_n}g}\inN$ is an unconditional frame for $M^p(\R)$ with the coefficient functionals $\set{g_n^*}\inN\subseteq M^{p'}(\R)$. Then we have 
$$f=\sumli\ip{f}{T_{\lambda_n}g}g_n^*$$
with unconditional convergence of the series in $M^{p'}(\R)$ for all $f\in M^{p'}(\R).$


\item [\textup{(d)}] (\cite[Theorem 11.3.2]{Gro01}) Fix $\psi\in S(\R)$ with $\norm{\psi}_{L^2(\R)}=1.$ Then for any $1\leq p<\infty$ we have $ \ip{f}{g} = \ip{V_\psi f}{V_\psi g}$ for any $f\in M^{p}(\R)$ and $g\in M^{p'}(\R).$

\item [\textup{(e)}] (\cite[Theorem 12.1.7]{Gro01}) For any $\phi \in M^1(\R)$ and $f\in M^p(\R)$ we have $$\norm{f\phi}_{M^p(\R)}\lesssim_p \norm{\phi}_{M^1(\R)}\norm{f}_{M^p(\R)}.$$

    \item [\textup{(f)}] (\cite[Proposition 1.7]{To04}) If $p\leq q\leq p'$, then $M^p(\R)\hookrightarrow L^q(\R)$. Consequently, we have that $M^1(\R)\hookrightarrow C_0(\R).$ \qeddef
    \end{enumerate}
    
\end{theorem}

For a fixed $y\in \R$ we define the corresponding modulation operator $M_y:M^p(\R)\rightarrow M^p(\R)$ by $(M_yg)(x)=e^{-2\pi iyx}g(x)$. It is known that for fixed $\alpha,\beta>0$ if $g,\gamma\in M^1(\R)$ are such that 
\begin{equation}
\label{Gabor_expansion}
f=\sum_{n,k\in \Z}\ip{f}{M_{\beta n}T_{\alpha k}\gamma}M_{\beta n}T_{\alpha m}g 
\end{equation}
with the convergence of the series in the norm of $M^2(\R)$ for all $f\in M^2(\R)$, then Equation (\ref{Gabor_expansion}) holds for all $f\in M^p(\R)$ with the unconditional convergence of the series for all $1\leq p<\infty.$ 
Yet, finding a pair $g$ and $\gamma$ in $M^1(\R)$ that satisfy Equation (\ref{Gabor_expansion}) simultaneously for all $f$ in $M^2(\R)$ is often challenging. There have been several known techniques for the construction of an unconditional frame for $M^2(\R)$ of the form $\set{M_{\beta n}T_{\alpha m}g}_{n,m\in \Z}$ in the field of \emph{Gabor analysis}. For example, the famous \emph{Painless Nonorthogonal Expansions} discovered by Daubechies, Grossman, and Meyer in \cite{DGM86} provides an useful tool to construct an unconditional frame of the form $\set{M_{\beta n}T_{\alpha m}g}_{n,m\in \Z}$ with $g$ being smooth and compactly supported. On the other hand, a deep result due to Gr\"{o}chenig and Leinart in \cite{GL03} showed that if $\set{M_{\beta n}T_{\alpha k}g}_{k,n\in \Z}$ with $g\in M^1(\R)$ is an unconditional frame for $M^2(\R)$ for some $\alpha,\beta>0$, then there exists at least one $\gamma\in M^1(\R)$ such that Equation (\ref{Gabor_expansion}) holds for all $f\in M^2(\R)$, and hence holds for all $f\in M^p(\R)$ with the unconditional convergence of the series for all $1\leq p<\infty.$ Combining these two results, we obtain the following theorem.  We will use $\norm{\cdot}_X\approx \norm{\cdot}_Y$ to mean that $\norm{\cdot}_X$ and $\norm{\cdot}_Y$ are equivalent. That is, there exist positive constants $A$ and $B$ such that $A\norm{\cdot}_Y\Le \norm{\cdot}_X\Le B\norm{\cdot}_Y$.

    \begin{theorem} (\cite[Theorem 1]{DGM86}, \cite[Corollary 12.2.6]{Gro01} and \cite[Theorem 4.2]{GL03}) 
	\label{painless_non_expan}
Fix any $0<\alpha<1$. There exists some $\phi\in C^\infty_c(\R)$ supported in $[0,1]$ such that for any $f\in M^p(\R)$ there exists a sequence of scalars $(c_{kn})_{k,n\in \Z}\in \ell^p(\Z^2)$ for which 
$$f=\sum_{k,n\in \Z}c_{kn}M_{n}T_{\alpha k}\phi$$
with the unconditional convergence of the series in $M^p(\R)$ for all $1\leq p<\infty.$  Moreover, we have $\norm{f}_{M^p(\R)}\approx \norm{(c_{kn})_{k,n\in \Z}}_{\ell^p(\Z^2)}\approx \norm{\bigparen{\ip{f}{M_nT_{\alpha k}\phi}}_{k,n\in \Z}}_{\ell^p(\Z^2)}$  \qeddef
\end{theorem}

In \cite{DH00}, Deng and Heil proved a necessary condition for a sequence of translates $\schtr$ to be an Schauder basis for a family of Banach spaces. 
Since modulation spaces fall into this category, we obtain the following special case of \cite[Lemma 4.1]{DH00}. 
 We say that a subset $\Lambda=\set{\lambda_n}\inN$ of $\R$ is uniformly discrete if $\min_{n\neq m} |\lambda_n-\lambda_m|>\delta$ for some $\delta>0.$
\begin{lemma}(\cite[Lemma 4.1]{DH00})
\label{Schauder_discre_index}
Let $1\leq p<\infty$. 
    Assume that $\schtr$ is a Schauder basis for $M^p(\R)$. Then $\set{\lambda_n}\inN$ is uniformly discrete. \qeddef
\end{lemma}

Finally, we state the boundedness of the discrete Hilbert transform. For a proof, we refer to \cite[Chapter 13]{FK09}).
\begin{lemma}
	\label{discrete_hilbert}
	For each $1<p<\infty$ there exists a constant $C_p>0$ such that for any sequence of scalars $(c_n)_{\inN} \in \ell^p(\Z)$ we have $$\Bignorm{\Bigparen{\sum_{\substack{n\,\in \,\Z\\ n\neq m}}c_{n} \cdot \frac{1}{m-n}}_{m\in\Z}}_{\ell^p(\Z)}\Le C_p \,\bignorm{(c_n)_{n\in \Z}}_{{\ell^p}(\Z)}.$$
\end{lemma}

\section{System of Translates in Modulation Spaces} 
 The modulation space $M^1(\R)$ has a significantly different nature from $M^p(\R)$ with $p>1$. For example, $M^1(\R)$ is continuously embedded into $C_0(\R)$, while other modulations spaces do not satisfy this property.
 This significant difference suggests dealing with $M^1(\R)$ and $M^p(\R)$ with $1<p\leq 2$ separately. Consequently, we will start with the case $1<p\leq 2$ and show that no systems of translates can be an unconditional basis or an unconditional frame for $M^p(\R)$ for any $1<p\leq 2.$ Of course, that a sequence is not an unconditional frame implies that it is not an unconditional basis. Nevertheless, the notions of unconditional bases and unconditional frames are often of independent interest in the literature, therefore, we will utilize two different approaches to disprove the existence of an unconditional basis and the existence of an unconditional frame for $M^p(\R)$ with $1<p\leq 2,$ respectively. 

We will need a series of lemmas. The first lemma follows from a slight adaptation of the proof of \cite[Proposition 3.1.3]{AK06}. 

\begin{lemma} 
\label{uncondi_cons_unsch}
Let $X$ be a separable Banach space. Assume that $\set{x_n}\inN$ is an unconditional frame for $X$ with coefficient functionals $\set{x_n^*}\inN$. Then there exists a constant $C>0$ for which $$\Bignorm{\sumli c_n\ip{f}{x_n^*}\,x_n}_{X}\Le C\norm{f}_{X}$$
   for any sequence of scalars $(c_n)_{\inN}$ with $\norm{(c_n)_{\inN}}_{\ell^\infty(\N)}\leq 1$ and all $f\in X.$
\end{lemma}

Due the fact that $M^p(\R)$ is invariant under the Fourier transform, we obtain a necessary condition for a sequence of translates to be a Schauder frame for $M^p(\R)$ for $1\leq p\leq 2$ as follows.
\begin{lemma}
\label{nonzero_ae_atom_schtr}
    Fix $1\leq p\leq 2$. Let $\Lambda=\set{\lambda_n}_{\inN} $ be a subset of $\R$ and let $g\in M^p(\R)$. Assume that $\schtr$ is a Schauder frame for $M^p(\R)$. Then $\widehat{g}(x)\neq 0$ for almost every $x\in \R$.
\end{lemma}
\begin{proof} Let $\set{g_n^*}\inN$ be the associated coefficient functionals with $\schtr$.
    Suppose to the contrary that $Z=\set{x \in \R\,|\, \widehat{g}(x)=0}$ has a positive measure. Let $K>0$ be large enough that $Z_K=[-K,K]\cap Z$ has a positive measure. Then we choose $\phi \in S(\R)$ for which $\widehat{\phi}=1$ on $[-K,K]$ and $\widehat{\phi}=0$ on $(-\infty,-K-1]\cup [K+1,\infty)$. By Theorem \ref{tf_lemmas} (b) and (f), we see that 
     $\widehat{\phi}(x) = \sumli \ip{\phi}{g_n^*}\,e^{-2\pi i \lambda_n x} \,\widehat{g}(x) $
    with the convergence of the series in the norm of $L^p(\R)$. However, this would imply $$0=\sumli \ip{\phi}{g_n^*} \int_\R e^{-2\pi i\lambda_n x} \widehat{g}(x) \chi_{Z_K} (x) \,dx = \int_{[-K,K]} \widehat{\phi}(x) \chi_{Z_K} (x) dx= |Z_K|>0,$$
    which is a contradiction. 
\end{proof}

When we have an unconditional frame  $\set{T_{\lambda_n}g}\inN$ for $M^p(\R)$ with coefficient functionals $\set{g_n^*}\inN,$ it is natural to ask whether the coefficients $\bigparen{\ip{f}{g_n^*}}\inN$ are bounded in some sense for all $f\in M^p(\R)$ since it is closely related the convergence of $\sumli\ip{f}{g_n^*}T_{\lambda_n}g$.
The next lemma shows that we do have some controls over the coefficients under some conditions. 

Recall that Khintchine's Inequality states that for any $1\leq p<\infty$ there exists positive constants $A_p,B_p$ such that for every $N\in \N$ and scalars $c_1,\dots,c_N$ 
$$A_p\norm{(c_n)_{n=1}^N}_{\ell^2(N)} \Le \Bignorm{\sum_{n=1}^N c_nR_n}_{L^p([0,1])}\Le B_p\norm{(c_n)_{n=1}^N}_{\ell^2(N)}.$$
Here $\set{R_n(x)}_{n\geq 0}$ denotes the Rademacher system, which is defined by $R_n(x)= \text{sign}(\sin(2^n \pi x))$ for $n\geq 0$.

\begin{lemma}
\label{coef_contro_Mp_unfr}
      Fix $1\leq p<\infty.$ Assume that $\schtr$ is an unconditional frame for $M^p(\R)$ with coefficient functionals $\set{g_n^*}\inN$. Assume that $\set{\lambda_n}\inN$ is uniformly discrete. Then the following statements hold:
      \begin{enumerate}
          \item [\textup{(a)}] There exists a constant $C>0$ depending on $g$ and $p$ such that $$\bignorm{\bigparen{\ip{f}{g_n^*}}\inN}_{\ell^p(\N)}\Le C \norm{f}_{M^p(\R)},$$ for all $f\in M^p(\R).$
           \item [\textup{(b)}] If $2<p<\infty$, then there exists a constant $K>0$ depending on $g$ and $p$ such that $$K\norm{f}_{M^{p'}(\R)}\Le \bignorm{\bigparen{\ip{f}{T_{\lambda_n}g}}\inN}_{\ell^{p'}(\N)},$$ for all $f\in M^{p'}(\R).$
         \end{enumerate}
\end{lemma}
\begin{proof}
  Fix $\psi \in S(\R)$ with $\norm{\psi}_{L^2(\R)}=1.$  For notational convenience, for $n\in \N$ we denote $V_\psi (T_{\lambda_n}g), V_\psi (g_n^*)$ by $\psi_n,\psi^*_n$, respectively.
  
(a)   It is not hard to verify that $S=V_\psi(M^p(\R))$ is a Banach space equipped with $L^p$-norm. By Theorem \ref{tf_lemmas} (d), we see that 
   $ V_\psi f 
    = 
     \sumli \,\ip{V_\psi f}{\psi_n^*} \, \psi_n,$ for all $f\in M^p(\R)$.
So, $\bigset{\psi_n}\inN$ is an unconditional frame for $S$ with coefficient functionals $\set{\psi_n^*}\inN$. Then by monotone convergence theorem and Lemma \ref{uncondi_cons_unsch} we compute \begin{equation*}
    \begin{split}
        \bignorm{\bignorm{\bigparen{\ip{V_\psi f}{\psi_n^*}\,\psi_n}}_{\ell^2(\N)}}_{L^p(\R^2)} &\Eq \lim_{N\rightarrow \infty} \Bignorm{\Bigparen{\sum_{n=1}^N \,\bigabs{\ip{V_\psi f}{\psi_n^*}\,\psi_n(\cdot)}^2}^{1/2}}_{L^p(\R^2)}\\
        &\Le \lim_{N\rightarrow \infty}A_p^{-1}\Bignorm{\Bignorm{\sum_{n=1}^N\ip{V_\psi f}{\psi_n^*}\,\psi_n(\cdot )R_n}_{L^p(\R^2)}}_{L^p([0,1])}\\
        &\Le KA_p^{-1} \norm{\norm{V_\psi f}_{L^p(\R^2)}}_{L^p([0,1])}\\
        &\Eq KA_p^{-1} \norm{f}_{M^p(\R)},
    \end{split}
\end{equation*}
for some positive constants $K$ and $k_p.$
Here the second inequality is a direction application of Khintchine's inequality. Next, assume that $\min_{i\neq j}|\lambda_i-\lambda_j|> 2\delta >0$ for some $\delta>0$ and assume that $\norm{V_{\psi}g|_{[k_0\delta,\,(k_0+1)\delta]\times \R}}_{L^p(\R)}>0$ for some $k_0\in \Z$. Here $V_{\psi}g|_{[k_0\delta,\,(k_0+1)\delta]\times \R}$ denotes the restriction of $V_{\psi}g$ to $[k_0\delta,\,(k_0+1)\delta]\times \R.$ Note that since $|\lambda_i-\lambda_j|> 2\delta$, for $i\neq j$ the intervals
$[\lambda_n+k_0\delta,\, \lambda_n+(k_0+1)\delta]$ and $[\lambda_m+k_0\delta,\, \lambda_m+(k_0+1)\delta]$ are disjoint. Then we compute \begin{equation*}
    \begin{split}
         \Bignorm{\bignorm{\bigparen{\ip{V_\psi f}{\psi_n^*}\,\psi_n}}_{\ell^2(\N)}}_{L^p(\R^2)}^p &\Eq \int_{\R^2} \bignorm{\bigparen{\ip{V_\psi f}{\psi_n^*}\,\psi_n(\cdot)}}_{\ell^2(\N)}^p\\
         &\Ge \sum_{m=1}^\infty \int_{[\lambda_m+k_0\delta,\, \lambda_m+(k_0+1)\delta] \times \R}  \bigabs{\ip{f}{g_m^*} V_{\psi} (x-\lambda_m,w)}^p\\
         &\Eq \sum_{m=1}^\infty |\ip{f}{g_m^*}|^p\int_{[k_0\delta,\, (k_0+1)\delta] \times \R} \bigabs{ V_{\psi} (x,w)}^p\\
    \end{split}
\end{equation*}
Thus, we conclude that $\Bignorm{\bigparen{\ip{f}{g_n^*}}}_{\ell^p(\N)}\Le C \norm{f}_{M^p(\R)}$ for some $C>0.$\medskip

(b) By Theorem \ref{tf_lemmas}(c), we have that $\set{\psi_n^*}\inN$ is an unconditional frame for $V_{\psi}(M^{p'}(\R))$ with coefficient functionals $\set{\psi_n}\inN$. Similarly, for each $N\in \N$ we compute
\begin{align*}
    \begin{split}
         \Bignorm{\sum_{n=1}^N \ip{V_\psi f}{\psi_n} \psi^*_n }_{L^{p'}(\R^2)} &\Eq  \biggnorm{\Bignorm{\sum_{n=1}^N \ip{V_\psi f}{\psi_n} \psi^*_n }_{L^{p'}(\R^2)}}_{L^{p'}([0,1])}\\
        &\Le C\,\biggnorm{\Bignorm{\sum_{n=1}^N \ip{V_\psi f}{\psi_n} R_n(t)\psi^*_n }_{L^{p'}(\R^2)}}_{L^{p'}([0,1])}\\
        &\Eq CB_p\,\Bignorm{\Bigparen{\sum_{n=1}^N \bigabs{\ip{V_\psi f}{\psi_n}\psi^*_n}^2}^{1/2} }_{L^{p'}(\R^2)}\\
    \end{split}
\end{align*}
for some positive constants $C$ and $B_p$. Here the second inequality is a direct application of Lemma \ref{uncondi_cons_unsch}.
 Note that we have $\sup_n\norm{\psi_n^*}_{L^{p'}(\R)}<\infty$ since $\set{\psi_n^*}\inN$ converges weak$^{*}$ to $0$.
 Statement (b) then follows from the following inequality  
\begin{align*}
    \begin{split}
    \Bignorm{\Bigparen{\sum_{n=1}^N \bigabs{\ip{V_\psi f}{\psi_n}\psi^*_n}^2}^{1/2} }_{L^{p'}(\R^2)} &\Le
 \Bignorm{\Bigparen{\sum_{n=1}^N \bigabs{\ip{V_\psi f}{\psi_n}\psi^*_n}^{p'}}^{1/p'} }_{L^{p'}(\R^2)}\\ &\Le \sup\limits_{n\in\N} \norm{\psi_n^*}_{L^{p'}(\R)} \Bigparen{\sum_{n=1}^N |\ip{V_\psi f}{\psi_n}|^{p'}}^{1/p'}. \qedhere
    \end{split}
\end{align*}
\end{proof}

In practice, it is usually very hard to directly use Definition 2.1 to compute the $M^p$-norm of an arbitrary tempered distribution explicitly. We will need to compute the $M^p$-norm of every function in the Rademacher system. To do so, we will appeal to the following theorem, which is a special case of \cite[Theorem 3.5]{PY24}. We remark that the original proof can be significantly simplified in our settings. To make this paper self-contained, we include the simplified proof here.

    \begin{theorem}
 		\label{conv_Gabor_ext} Fix $0<\alpha\leq 1.$ Then for any $1<p<\infty$ the following statements hold:    
    \begin{enumerate}
       \setlength\itemsep{0.5em}
    \item [\textup{(a)}] The Gabor expansions \begin{equation}
      \label{Gabor_exp}
    \begin{split}
        f&\Eq\sum_{n,k\in \Z} \bigip{f}{M_{n}T_{\alpha k}\Chi}  M_{n}T_{\alpha k}\Chi
    \end{split}
     \end{equation}
     hold with unconditional convergence of the series in $M^{p}(\R)$ for all $f\in M^{p}(\R)$.
     \item [\textup{(b)}]
 			    There exist positive constants $A$ and $B$ depending on $p$ and $\alpha$ such that  $$A\,\norm{f}_{M^p(\R)}\Le\Bignorm{\Bigparen{\bigip{f}{M_{ n}T_{\alpha k}(\Chi)}}_{k,n\in \Z}}_{\ell^{p}(\Z^2)} \Le B\,\norm{f}_{M^{p}(\R)},$$
                for all $f\in M^p(\R).$
         \end{enumerate}
\end{theorem}
\begin{proof} 
    We first define the analysis operator $C_{\Chi}\colon M^p(\R) \rightarrow \ell^p(\Z^2)$ associated with $\Chi$ by $$C_{\Chi}(f)\Eq\bigparen{\ip{f}{M_{n}T_{\alpha k} g}}_{k,n\in \Z}$$
    and define the synthesis operator associated with $R_{\Chi}\colon\ell^p(\Z^2)\rightarrow M^p(\R)$ by $$R_{\Chi}(c_{kn})= \sum_{n,k\in \Z} c_{kn}  M_{n}T_{\alpha k}\Chi.$$ 
We will show that $C_{\Chi}\colon M^p(\R) \rightarrow \ell^p(\Z^2)$ is a bounded linear operator for all $1<p<\infty.$ 
The boundedness of the synthesis operator $R_{\Chi}\colon\ell^p(\Z^2)\rightarrow M^p(\R)$ for all $1<p<\infty$ then follows from the fact that $R_{\Chi}$ is the adjoint operator of $ C_{\Chi}.$ We will then utilize the boundedness of $C_{\Chi}$ and $R_{\Chi}$ combined with the fact that $\set{M_{n}T_{\alpha k}\Chi}_{n,k\in \Z}$ is an unconditional frame for $M^2(\R)$ with coefficient functionals $\set{M_{n}T_{\alpha k}\Chi}_{n,k\in \Z}$ to extend Equation (\ref{Gabor_exp}) to all $1<p<\infty$. Finally, statement (b) follows by the following inequality
$$\norm{f}_{M^{p}(\R)} \Eq \norm{(R_{\Chi}\circ C_{\Chi})(f)}_{M^{p}(\R)} \Le  K\norm{C_{\Chi}(f)}_{\ell^{p,q}(\Z^2)}   \Le K' \norm{f}_{M^{p,q}(\R)},$$
for some constants $K$ and $K'$ depending only on $p$ and $\alpha$. 

Therefore, it remains to show that $C_{\Chi}:M^p(\R) \rightarrow \ell^p(\Z^2)$ is a bounded linear operator for any $1<p<\infty$.
By Theorem \ref{painless_non_expan}, there exists some $\phi \in C^\infty_c(\R)$ supported in $[0,\alpha]$ such that every $f\in M^p(\R)$ can be expressed as 
$$f\Eq \sumli c_{mj} M_{j}T_{\frac{\alpha m}{2}}\phi$$ with unconditional convergence of the series in $M^p(\R)$. Here $(c_{mj})_{m,j\in \Z}$ is a sequence of scalars in $\ell^p(\Z^2)$ for which $\norm{(c_{mj})_{m,j\in \Z}}_{\ell^p(\Z^2)}\approx \norm{f}_{M^p(\R)}$. Since $\Chi\in M^{p}(\R)$ for any $p>1$, we see that for any fixed $k,n\in \Z$ the sequence $\bigparen{\bigip{ M_{m} T_{\frac{\alpha j}{N}} \phi}{M_{k}T_{\alpha n}\Chi}}_{j,m\in \Z}$ is in $\ell^p(\Z^2)$ for any $p>1$ by Theorem \ref{painless_non_expan}. Using H\"{o}lder's inequality to justify the interchange of inner product and summation, we obtain 
$$\bigabs{\ip{f}{M_{ n}T_{\alpha k}\Chi} }
	\Eq \Bigabs{\sum_{m,j\in \Z} c_{mj}\bigip{ M_{j} T_{\frac{\alpha m}{2}} \phi}{M_{n}T_{\alpha k}\Chi}}.$$
Note that $\bigip{ M_{j} T_{\frac{\alpha m}{2}} \phi}{M_{n}T_{\alpha k}\Chi}$ is nonzero only when $m=2k,2k-1$ and $2k+1$. For the case $j=n$, we have$$\Bignorm{\bigparen{\ip{f}{M_{n}T_{\alpha k} g}}_{k,n\in \Z}}^p_{\ell^p(\Z^2)}\Le C \sum_{k,n\in \Z}\sum_{2k-1\leq m\leq 2k+1} |c_{mn}|^p\Le C'\,\norm{(c_{kn})}_{\ell^p(\Z^2)}^p$$
for some constant $C$ and $C'>0$ by the Triangle Inequality. 

Next, we address the case $j\neq n.$ We first estimate the case $m=2k$. For each $j,n\in \Z$, we let $K_{j-n} = e^{ i 2\pi (j-n)(k+1)\alpha} \phi (\alpha )$, $L_{j-n} = e^{ i 2\pi (j-n)k\alpha} \phi (0)$ and $$H_{j-n}=\int^{\alpha (k+1)}_{\alpha k} e^{ 2\pi (j-n)x} \phi' (x-\alpha k)\, dx.$$By integration by parts, we see that 
\begin{align*}
\begin{split}
\bigip{ M_{j} T_{\frac{\alpha m}{2}} \phi}{M_{n}T_{\alpha k}\Chi} &\Eq \int^{\alpha (k+1)}_{\alpha k} e^{ 2\pi i (j-n)x}  \phi(x-\alpha k) \,dx\\
&\Eq \frac{K_{j-n}-L_{j-n}+H_{j-n}}{ 2\pi i (j-n)}.
\end{split}    
\end{align*}
Utilizing the boundedness of discrete Hilbert transform, we have 
$$ \Bignorm{\Bigparen{\sum_{j\in \Z} c_{(2k)j}\dfrac{K_{j-n}}{(j-n)}}_{k,n\in \Z}}_{\ell^p(\Z^2)} \lesssim_{\phi,p,\alpha} \norm{(c_{mj})_{m,j\in \Z}}_{\ell^p(\Z^2)},$$
and 

$$ \Bignorm{\Bigparen{\sum_{j\in \Z} c_{(2k)j}\dfrac{L_{j-n}}{(j-n)}}_{k,n\in \Z}}_{\ell^p(\Z^2)} \lesssim_{\phi,p,\alpha} \norm{(c_{mj})_{m,j\in \Z}}_{\ell^p(\Z^2)}.$$
Using integration by parts again, we see that $\bigabs{\frac{H_{j-n}}{j-n}}$ is dominated by $\frac{1}{(j-n)^2}$ multiplied by some constant depending on $\phi, p$ and $\alpha.$ Therefore, 
$$ \Bignorm{\Bigparen{\sum_{j\in \Z} c_{(2k)j}\dfrac{H_{j-n}}{(j-n)}}_{k,n\in \Z}}_{\ell^p(\Z^2)} \lesssim_{\phi,p,\alpha} \norm{(c_{mj})_{m,j\in \Z}}_{\ell^p(\Z^2)},$$
by Young's convolution inequality. Using the Triangle Inequality, we reach the conclusion that $$\bignorm{C_{\Chi}(f)}\lesssim_{\phi,p,\alpha} \norm{f}_{M^{p}(\R)}. $$ 
The case $m=2k-1$ and $m=2k+1$ can be proved similarly.
\qedhere   
\end{proof}

\begin{lemma}
\label{uniform_bound_Mpnorm_Ra}
Let $\set{R_n}_{n\in \N\cup \set{0}}$ be the Rademacher system. Then for any $1<p<\infty$ there exist positive constants $A_p$ and $B_p$ depending on $p$ only such that  $$A_p\norm{R_n}_{M^p(\R)} \Le  \bignorm{\bigparen{\frac{1}{j}}_{j \in \N}}_{\ell^p(\N)} \Le B_p\norm{R_n}_{M^p(\R)}$$
for all $n\in \N.$
\end{lemma}
\begin{proof} Fix $N\in\N.$
    By Theorem \ref{conv_Gabor_ext}, there exist positive constants $A_p$ and $B_p$ such that 
  \begin{equation}
  \label{norm_equiv_Rade}
            A_p \norm{R_N}_{M^p(\R)} \Le  \Bignorm{\Bigparen{\bigip{R_N}{M_{ k}T_{\ell}(\chi_{_{[0,1]}})}}_{k,n\in \Z}}_{\ell^{p}(\Z^2)} \Le B_p\norm{R_N}_{M^p(\R)}
 \end{equation}
Note that $\bignorm{\bigparen{\bigip{R_N}{M_{k}T_{\ell}(\chi_{_{[0,1]}})}}_{k,\ell\in \Z}}_{\ell^{p}(\Z^2)}= \bignorm{\bigparen{\bigip{R_N}{M_{k}(\chi_{_{[0,1]}})}}_{k\in \Z}}_{\ell^{p}(\Z)}$. Then for $-k\neq 0$ we compute \begin{align*}
        \begin{split}
        \ip{R_N}{e_{-k}} &\Eq \sum_{n=0}^{2^N-1} (-1)^n\int^{\frac{n+1}{2^N}}_{\frac{n}{2^N}}e^{2\pi ikt} dt \\
        &\Eq \sum_{n=0}^{2^N-1}(-1)^n \frac{1}{2\pi i k}[e^{2\pi ik(\frac{n+1}{2^N})}-e^{2\pi ik(\frac{n}{2^N})}]\\
        &\Eq \frac{1}{\pi ik}\sum_{n=0}^{2^N-1}(-1)^{n+1} e^{2\pi i \frac{n}{2^N}k}\\
         \end{split}
    \end{align*}
If $k\notin \set{j2^{N-1} \,|\,j\in \Z \text{ is an odd number}}$, then we have $\frac{1}{\pi ik}\sum_{n=0}^{2^N-1}(-1)^{n+1} e^{2\pi i \frac{n}{2^N}k} =0$.
If $k=j2^{N-1}$ for some odd integer $j$, then we have 
\begin{align*}
    \begin{split}
   \frac{1}{\pi ik}\sum_{n=0}^{2^N-1}(-1)^{n+1} e^{2\pi i \frac{n}{2^N}k} &\Eq \dfrac{1}{\pi i 2^{N-1}\ell} \sum_{n=0}^{2^N-1}(-1)^{n+1} e^{inj\pi}\\
   &\Eq \dfrac{1}{\pi i 2^{N-1}j} \sum_{n=0}^{2^N-1}(-1)^{2n+1}\\
   & \Eq \dfrac{-2}{\pi ij}.
\end{split}
   \end{align*}
The statement then follows from Equation (\ref{norm_equiv_Rade}).  
\end{proof}

\begin{corollary}
\label{decay_lemma_Rade}
    Let $\set{R_n}_{n\in \N\cup \set{0}}$ be the Rademacher system. Then for any $1<p<\infty$ we have $\lim\limits_{n\rightarrow \infty} \ip{R_n}{f}=0$  for all $f\in M^p(\R).$
\end{corollary}
\begin{proof}
    Note that $\ip{R_n}{f}$ is well-defined for all $n\in \N\cup \set{0}$ since $R_n \in M^p(\R)$ for any $1<p<\infty$. Let $\set{f_k}_{k\in \N}\in S(\R)$ be a sequence that converges to $f$ in $M^p(\R)$. It is well-known that $\ip{R_n}{f_k}\rightarrow 0$ as $n\rightarrow \infty$ for any $k\in \N$. By the Triangle Ineqaulity, we obtain $$|\ip{R_n}{f}|\Le \bigparen{\sup_{n \in \N}\norm{R_n}_{M^{p'}(\R)}}\norm{f_k-f}_{M^{p}(\R)}+|\ip{R_n}{f_k}|.$$
    The result then follows from Lemma \ref{uniform_bound_Mpnorm_Ra}.
\end{proof}

For any $a<b$, we define $C^1_{AC}[a,b]$ to be the set of all complex-valued supported in $[a,b]$ that have absolutely continuous derivatives on $[a,b]$, i.e.,
$$C^1_{AC}[a,b] = \bigset{g\colon \R\rightarrow \C\,\big|\,\textup{supp\,}(g)\subseteq [a,b] \text{ and }\,g'(x) \text{ is absolutely continuous on }[a,b] }.$$
We remark that $C^1_{AC}[a,b]\subsetneq M^p(\R)$ for any $p>1$ by Theorem \ref{conv_Gabor_ext}.
\begin{lemma}
\label{discrete_infinite_sum_Mpnorm} 
Fix $a>b$ and fix $\psi\in C^1_{AC}[a,b]\cup M^1(\R)$. Assume that $\Lambda=\set{\lambda_n}\inN \subseteq \R$ is a discrete subset. Then for any $1<p<\infty$ we have 
$$\sumli \norm{T_{\lambda_n}f\cdot \psi\cdot \chi_{_{[a,b]}}}_{M^p(\R)}^p<\infty,$$
for all $f\in M^p(\R).$
\end{lemma}
\begin{proof}
We prove the case that $\psi\in C^1_{AC}[a,b]$ first.
Let $\delta>0$ be such that $\min\limits_{i\neq j} |\lambda_i-\lambda_j|>2\delta.$ Using the Triangle Inequality if necessary, we may assume $b-a<2\min \set{\delta,1}.$ By Theorem \ref{conv_Gabor_ext}, we have 
\begin{align*}
\begin{split}
    \sumli \norm{T_{\lambda_n}f\cdot \psi }_{M^p(\R)}^p &\approx \sumli \Bigparen{\sum_{k,\ell \in \Z}^\infty \bigabs{\ip{T_{\lambda_n}f \cdot \psi}{M_kT_{(b-a)\ell}(\chi_{_{[a,b]}})}}^p} \\
 &\Eq \sumli \sum_{k=1}^\infty\, \bigabs{\ip{f }{M_kT_{-\lambda_n}\overline{\psi}}}^p.
    \end{split}
\end{align*}
Using Thereom \ref{conv_Gabor_ext} again and use H\"{o}lder's Ineqaulity to justify the interchange of summation and inner product, we see that 
\begin{align*}
    \ip{f }{M_kT_{-\lambda_n}\overline{\psi}} \Eq \sum_{m,j\in \Z} c_{mj} \ip{M_jT_{(b-a)m}\chi_{_{[a,b]}}}{M_kT_{-\lambda_n}\overline{\psi}}.
\end{align*} 
Here $(c_{m_j})_{m,j\in \Z}$ is a sequence of scalars for which $\norm{(c_{m_j})_{m,j\in \Z}}_{\ell^p(\Z^2)}\approx \norm{f}_{M^p(\R)}.$
Fix $k$ and $j\in \Z$. For each $n\in\Z$ we define $$S^{k,j}_n = \set{m\in \Z\,|\, \ip{M_jT_{(b-a)m}\chi_{_{[a,b]}}}{M_kT_{-\lambda_n}\overline{\psi}} \neq 0}.$$
Note that the cardinality of $S^{k,j}_n$ is independent of $k$ and $j$, i.e., $|S^{k_1,j_1}_n|=|S^{k_2,j_2}_n|$ for any $k_1,j_1,k_2,j_2\in \Z$. Consequently, for notational convenience, we will simply write $S_n$ instead of $S^{k,j}_n.$ Since $\chi_{_{[a,b]}}$ is compactly supported, we see that $\sup_n |S_n|<\infty.$ Using the assumption that $b-a<2\delta$, we also have $S_{n_1}\cap S_{n_2}= \emptyset$ if $n_1\neq n_2$, and hence $\cup_{n\in \Z} S_n \subseteq \Z.$ Next, by the Triangle Inequality, we obtain
\begin{align*}
    \begin{split}
 \bigabs{\ip{f }{M_kT_{-\lambda_n}\overline{\psi}}} \Le 2^{p\sup_{n}|S_n|} \sum_{m\in S_n}\Bigabs{\sum_{j\in \Z} c_{mj} \ip{M_{j-k}T_{(b-a)m+\lambda_n}\chi_{_{[a,b]}}}{\overline{\psi}}}^p.
\end{split}
\end{align*}
Similar to the estimate used in Theorem \ref{conv_Gabor_ext}, we first write  \begin{align*}
\begin{split}
\ip{M_{j-k}T_{(b-a)m+\lambda_n}\chi_{_{[a,b]}}}{\overline{\psi}} \Eq \int^{b+(b-a)m+\lambda_n}_{a+(b-a)m+\lambda_n} e^{2\pi i(j-k)x}\psi\,dx.
\end{split}
\end{align*}
Using integration by parts and argue similarly to proof of Theorem \ref{conv_Gabor_ext}, we obtain 
\begin{align*}
    \begin{split}
    \sum_{k\in \Z}\Bigabs{\sum_{j\in \Z} c_{mj} \ip{M_{j-k}T_{(b-a)m+\lambda_n}\chi_{_{[a,b]}}}{\overline{\psi}}}^p\lesssim_{\psi,p,\Lambda} \norm{(c_{m,\cdot})_{m\in \Z}}_{\ell^p(\Z)}^p.
\end{split}
\end{align*}
Finally, since $\cup_{n\in \Z} S_n \subseteq \Z$, we reach the estimate   
\begin{align*}
\begin{split}
\sumli \norm{T_{\lambda_n}f\cdot \psi}_{M^p(\R)}^p &\Les_{\psi,p,\Lambda} \sum_{\substack{m\in S_n\\ n\in \Z}} \norm{(c_{m,\cdot})_{m\in \Z}}_{\ell^p(\Z)}^p\\
&\Les_{\psi,p,\Lambda}  \norm{(c_{m,j})_{m,j\in \Z}}_{\ell^p(\Z^2)}\\
&\Les_{\psi,p,\Lambda} \norm{f}_{M^{p}(\R)}^p.
\end{split}
\end{align*}
Now if $\psi \in M^1(\R)$, then by Theorem \ref{tf_lemmas} (e) we have $$\sumli \norm{T_{\lambda_n}f\cdot \psi \cdot \chi_{_{[a,b]}}}_{M^p(\R)}^p\Le \Bigparen{\sumli \norm{T_{\lambda_n}f\cdot \chi_{_{[a,b]}}}_{M^p(\R)}^p}\norm{\psi}_{M^1(\R)}^p.$$
The statement then follows from the fact that $\chi_{_{[a,b]}}\in C^1_{AC}[a,b].$\qedhere
\end{proof}
\begin{remark}
We can see a significant difference between $L^p$ spaces and modulation spaces simply from the proof this lemma. The $L_p$-version of lemma \ref{discrete_infinite_sum_Mpnorm} was first proved in \cite[Proposition 2.1]{OSSZ11} by using the fact that $L_p$ spaces are solid Banach spaces for all $1\leq p \leq \infty$. 
In general, modulation spaces are not solid Banach spaces except $M^2(\R)$. As a result, 
it is much more uncertain whether any restriction of a tempered distribution in $M^p(\R)$ to a bounded interval is still in the same modulation space. Lemma \ref{discrete_infinite_sum_Mpnorm} confirms that this the case if $1<p<\infty$. Moreover, we emphasize that this lemma can cannot be extended to $M^1(\R)$ since $M^1(\R)$ is continuously embedded into $C_0(\R)$. \qeddef
\end{remark}

We are now ready to disprove the existence of unconditional Schauder bases formed by a system of translates in certain modulation spaces.

\begin{theorem}
\label{non_exist_discrete_unschtr}
Let $\Lambda=\set{\lambda_n}\inN$ be an arbitrary subset of $\R$. Then the following statements are true.
\begin{enumerate} \setlength\itemsep{0.5em}
    \item [\textup{(a)}] For any $g\in M^p(\R)$ with $1\leq p\leq 2$, the sequence of translates $\set{T_{\lambda_n}g}\inN$ is never an unconditional basis for $M^q(\R)$ for any $p\leq q\leq p'.$
    \item [\textup{(b)}] For any $g\in M^p(\R) \cap L^{q}(\R)$ with $2\leq p<\infty$ and $1\leq q\leq p'$, the sequence of translates $\set{T_{\lambda_n}g}\inN$ is never an unconditional basis for $M^p(\R)$.
    \end{enumerate}
\end{theorem}
\begin{proof}
(a)  Suppose to the contrary that $\set{T_{\lambda_n}g}\inN$ is an unconditional basis for $M^q(\R)$ with the coefficient functionals $\set{g_n^*}\inN$. Then $\Lambda$ is uniformly discrete by Lemma \ref{Schauder_discre_index}. By Lemma \ref{discrete_infinite_sum_Mpnorm} and Lemma \ref{uniform_bound_Mpnorm_Ra}, there exists $K$ large enough that \begin{equation}
        \label{Eq1}
    \Bigparen{\sum_{n=K}^\infty \norm{T_{\lambda_n}g \cdot \chi_{_{[0,1]}}}_{M^{q'}(\R)}^{q'}}^{1/q'}\sup\limits_{n\in \N} \norm{R_n}_{M^q(\R)}< \frac{1}{2}\inf\limits_{n\in \N} \norm{R_n}_{M^q(\R)}.
    \end{equation}
    Also, by Lemma \ref{decay_lemma_Rade}, there exists $N_0$ large enough such that 
    \begin{equation}
    \label{Eq2}
    \sup_{n\in \N}\norm{T_{\lambda_n}g\cdot \chi_{_{[0,1]}}}_{M^q(\R)} \sum_{n=1}^{K-1}|\ip{R_{N_0}}{g_n^*}|<\frac{1}{2}\inf\limits_{n\in \N} \norm{R_n}_{M^q(\R)}.
    \end{equation}
    Combining Equation (\ref{Eq1}), Equation (\ref{Eq2}) with the Triangle Inequality and H\"{o}lder's inequality we obtain 
\begin{align*}
    \begin{split}
        \norm{R_{N_0}}_{M^q(\R)}  &= \Bignorm{\sumli \ip{R_{N_0}}{g_n^*}T_{\lambda_n}g\cdot \chi_{_{[0,1]}}}_{M^q(\R)} \\
        &\Le \Bignorm{\sum_{n=1}^{K-1} \ip{R_{N_0}}{g_n^*}T_{\lambda_n}g\cdot \chi_{_{[0,1]}}}_{M^q(\R)} + \sum_{n=K}^\infty \bignorm{\ip{R_{N_0}}{g_n^*}T_{\lambda_n}g\cdot \chi_{_{[0,1]}}}_{M^q(\R)}\\
        &<\inf_{n\in\N} \norm{R_n}_{M^q(\R)},
    \end{split}
\end{align*}
which is a contradiction. 
\medskip

(b) Assume that $\set{g(x-\lambda_n)}\inN$ is an unconditional basis for $M^p(\R)$. 
Then by Lemma \ref{coef_contro_Mp_unfr} (b) and H\"{o}lder's inequality, there exists some constant $A>0$ such that for any $k,L\in \N$
\begin{align*}
\begin{split}
   A \norm{R_k}_{M^{p'}(\R)} &\Le \Bigparen{\sum_{n=1}^{L-1}\bigabs{\ip{R_k}{T_{\lambda_n}g}}^{p'}}^{1/p'}+ \Bigparen{\sum_{n=L}^\infty \bigabs{\ip{R_k}{T_{\lambda_n}g}}^{p'}}^{1/p'}\\
     &\Le \Bigparen{\sum_{n=1}^{L-1}\bigabs{\ip{R_k}{T_{\lambda_n}g}}^{p'}}^{1/p'}+ \Bigparen{\sum_{n=L}^\infty \bigabs{\ip{R_k}{T_{\lambda_n}g}}^{q}}^{1/q}\\
   &\Le \Bigparen{\sum_{n=1}^{L-1}\bigabs{\ip{R_k}{T_{\lambda_n}g}}^{p'}}^{1/p'}+\Bigparen{\sum_{n=L}^{\infty}\norm{T_{\lambda_n}g\cdot \chi_{_{[0,1]}}}^{q}_{L^{q}(\R)}}^{1/q}
    \end{split}
\end{align*}
Choosing $L$ large enough and using Lemma \ref{uniform_bound_Mpnorm_Ra}, we obtain a contradiction.
\end{proof}

Although we only stated the results for unconditional bases in Theorem \ref{non_exist_discrete_unschtr}, the proof of Theorem \ref{non_exist_discrete_unschtr} applies to unconditional frames formed by discrete translates as well.
Next, we show that no system of translates can be an unconditional frame for $M^p(\R)$ for any $1<p\leq 2.$ The proof of the following theorem was inspired by \cite[Theorem 10.28]{Hei11}, where the author credited A. Olevskii for the proof. For any real numbers $a<b$ we use $M^1([a,b])$ to denote the subset of $M^1(\R)$ that contains all tempered distributions that are supported in $[a,b]$.

\begin{theorem}
\label{vanish_positive_measure_thm}
   Fix $1<p\leq 2.$ Assume that $\set{T_{\lambda_n}g}\inN$ is an unconditional frame for its closed span. Assume further that for some $a<b$ there exists some $\psi\in C^1_{AC}[a,b]\cup M^1([a,b])$ for which $g\ast \widehat{\psi}\neq 0\in \clspan{(T_{\lambda_n} g)}\inN$. Then $\widehat{g}$ must vanish in a set of positive measure, i.e., $|\set{x\in \R\,|\,\widehat{g}(x)\neq 0}|>0$.  

   Consequently, $M^p(\R)$ does not admit any unconditional frame formed by a sequence of translates for any $1< p\leq 2.$
\end{theorem}
\begin{proof} Suppose to the contrary that $\widehat{g}$ is nonzero almost everywhere.
    Let $G= \widehat{g}\cdot \psi$. Note that $G$ is still in $M^p(\R)$ by Lemma $\ref{discrete_infinite_sum_Mpnorm}.$
By assumptions, there exists some sequence of scalars $(c_n)_{\inN}$ such that 
    \begin{equation}
    \label{rec_one}
    G = \sumli c_ne^{-2\pi i\lambda_nx}\widehat{g}(x)
    \end{equation}
    with the unconditional convergence of the series in $M^p(\R).$ For each $k\in \N$ let \begin{equation}
    \label{rec_two}
    G_k\Eq \sumli c_n e^{-2\pi i\lambda_n(x-k)}\widehat{g}(x).
    \end{equation}
    Then by Lemma \ref{uncondi_cons_unsch} there exists some constant $C>0$ for which $\norm{G}_{M^p(\R)}\leq C \norm{G_k}_{M^p(\R)}$ for all $k\in \N.$  
     
     Next, let $S_N(x)= \sum_{n=1}^N  c_n\,e^{-2\pi i\lambda_nx}$. Using the fact that $M^p(\R)\hookrightarrow L^p(\R)$, we obtain a subsequence $\set{N_j}_{j\in\N}$ for which 
    $S_{N_j}(x)\widehat{g}(x)$ and $S_{N_j}(x-k)\widehat{g}(x)$ converges to $G$ and $G_k$ pointwise almost everywhere, respectively. Since $\widehat{g}$ is nonzero almost everywhere, the series $S_{N_j}$ converges to $\psi$ pointwise almost everywhere. Therefore, $S_{N_j}(x-k)\widehat{g}$ converges pointwise to $\widehat{g} \cdot T_{k}\psi$ almost everywhere, and hence $G_k=\widehat{g} \cdot T_{k}\psi $ almost everywhere. Since $\norm{G}_{M^p(\R)}\leq C\norm{G_k}_{M^p(\R)}$ for all $k\in \N$, it remains to show that $\norm{\widehat{g} \cdot T_{k}\psi}_{M^p(\R)}$ tends to $0$ as $k\rightarrow \infty$ to obtain a contradiction. By Lemma \ref{discrete_infinite_sum_Mpnorm}, we have that 
    $$\sum_{n\in \Z} \norm{\widehat{g}\cdot T_{k}\psi }_{M^p(\R)}^p =   \sumli \norm{T_{-k}\widehat{g}\cdot \psi}_{M^p(\R)}^p<\infty,$$
which implies that $\lim_{k\rightarrow \infty} \norm{G_k}_{M^p(\R)}=0$, which is a contradiction. Consequently, we conclude that $M^p(\R)$ does not admit any unconditional frames formed by a system of translates for any $1<p\leq 2$ by Lemma \ref{nonzero_ae_atom_schtr}.
\end{proof}

\begin{remark}
Recall that we say a Banach space $X$ is a \emph{left Banach A-module} with respect to multiplication `$\cdot$' if $A$ is a Banach algebra and $\cdot:A\times X\rightarrow X$ is such that $\norm{a\cdot x}_X\leq \norm{a}_A\norm{x}_{X}$ for all $a\in A$ and $x\in X.$ By Theorem \ref{tf_lemmas} (b) and (e), we see that $M^1(\R)$ is a Banach algebra with respect to convolution as well as with respect to pointwise multiplication.
 Now, let $S$ be a closed subspace of $M^p(\R)$ that for some $1<p\leq 2.$ Assume further that $S$ admits an unconditional frame formed by a system of translates. Then
Theorem \ref{vanish_positive_measure_thm} provides an useful criterion to check whether $S$ is a Banach $M^1(\R)$-module with respect to convolution. Furthermore, assume that $E$ is a closed subspace that admits an unconditional frame formed by a weighted exponential system $\set{e^{2\pi i\lambda_n x}g(x)}\inN$ for some $g\in M^p(\R)$ for some $1<p\leq 2$ and $\set{\lambda_n}\inN\subseteq \R.$ Then Theorem \ref{vanish_positive_measure_thm} implies that $E$ is not a Banach $M^1(\R)$-module with respect to pointwise multiplication if the Fourier transform of $g$ is compactly supported. 
\begin{example}
Let $\phi=e^{-\pi x^2}$ be the Gaussian function. Then for any $\set{a_n}\inN$ such that $\set{T_{a_n}\phi}\inN$ is an unconditional frame for its closed span we have that $\clspan{\set{T_{a_n}\phi}}\inN$ is not a Banach $M^1(\R)$-module with respect to convolution. For example, let $A$ and $B$ be two positive constants such that $A\Le \bigparen{\sumli |\phi(x+n)|^2}^{1/2} \Le B$ for every $x\in [0,1]$. Then for any $\set{a_n}\inN\subseteq \R$ for which $\Bignorm{\Bigparen{\norm{T_n\phi-T_{a_n}\phi}_{L^2(\R)}}_{n\in \Z}}_{\ell^2(\Z)}<A$ we have that $\set{T_{a_n}\phi}\inN$ is an unconditional frame for its closed span (for example, see \cite[Section 10.4]{Hei11} and \cite[Theorem 22.2.1]{Chr16}). Using Theorem \ref{vanish_positive_measure_thm}, we see that $\clspan{\set{T_{a_n}\phi}}\inN$ is not a Banach $M^1(\R)$-module with respect to convolution.
\end{example}

\begin{example} Fix $1<p\leq 2$ and let $g \in M^p(\R)$ be such that $\Fc^{-1}g$ is compactly supported. Here $\Fc^{-1}g$ denotes the inverse Fourier transform of $g$. Assume that there exists a subset $\set{\lambda_n}\inN$ of $\R$ such that $\set{e^{2\pi i\lambda_n x}g(x)}$ is an unconditional frame for its closed span. Since $M^p(\R)$ is invariant under Fourier transform, we see that $\bigset{T_{-\lambda_n}(\Fc^{-1}g)}\inN$ is an unconditional frame for its closed span in $M^p(\R)$.  Since $g$ is nonzero almost everywhere, we see that $\clspan{\bigset{T_{-\lambda_n}(\Fc^{-1}g)}\inN}$ is not a Banach $M^1(\R)$-module with respect to convolution. Consequently, we obtain that $\clspan{\set{e^{2\pi i\lambda_n x}g(x)}}\inN$ is not a Banach $M^1(\R)$-module with respect to multiplication.
    
\end{example}
\end{remark}

Next, we study the existence of unconditional basis formed by a sequence of translates in $M^1(\R).$ Recall that 
the \emph{effective density}, $\tilde{D}_\Lambda$, associated with a discrete subset $\Lambda\subseteq \R$ is defined as follows. Let $\mathcal{U}=\set{(a_n,b_n]}_{\inN}$ be a sequence of non-overlapping increasing intervals contained in $[0,\infty)$, i.e., $0\leq a_0<b_0\leq a_1<b_1\leq a_2<b_2\cdots.$ Following the terminology in \cite[Chapter IX]{Ko92}, we say that $\mathcal{U}$ is \emph{substantial} for a positive number $C$ (with respect to $\Gamma $) if $$\frac{|(a_n,b_n]\cap \Gamma  |}{b_n-a_n}> C, \quad \text{and} \quad \sumli \,\bigparen{\frac{b_n-a_n}{a_n}}^2=\infty,$$
where $|(a_n,b_n]\cap \Gamma |$ denotes the number of elements of $\Gamma$ that are contained in $(a_n,b_n].$
The effective density $\tilde{D}_\Lambda$ is then defined as  the supremum of all numbers $C$ for which $\mathcal{U}$ is substantial. For general subset $\Lambda$ of real numbers, the effective density of $\Lambda$ is defined to be$$\tilde{D}_\Lambda = \max (\tilde{D}_{\Lambda^+}, \tilde{D}_{\Lambda^-}),$$
where $\Lambda^+=\R^+\cap \Lambda$ and $\Lambda^-=\R^+\cap (-\Lambda)$. The Beurling--Malliavin theorem of radius of completeness (\cite{BM67}, see also \cite{Ko92} for a detailed exposition) states that $$\sup\set{R>0\,|\,\set{e^{ i\lambda_nx }}_{\inN}~\text{is complete in }C([-R,R])}\Eq \pi \tilde{D}_\Lambda.$$

\begin{proposition}
\label{effective_density_fr_M1}
Let $\Lambda=\set{\lambda_n}\inN$ be an arbitrary subset of $\R.$ Assume that $\set{T_{\lambda_n}g}\inN$ is an unconditional frame for $M^1(\R).$ Then we have $$\sup\set{R>0\,|\,\set{e^{ -2\pi i\lambda_nx }}_{\inN}~\text{is complete in }C([-R,R])}\Eq \infty.$$
\end{proposition}

\begin{proof}
Let $c\in \R$ be such that $\widehat{g}(c)\neq 0.$ Assume that there exists some $\set{g_n^*}\inN\subseteq (M^1(\R))^*$ such that for all $f\in M^1(\R)$ we have \begin{equation}
\label{eq_uncon_fr_M1}
    f=\sumli \ip{f}{g_n^*}T_{\lambda_n}g
\end{equation}
with the unconditional convergence of the series in $M^1(\R)$. By Lemma \ref{uncondi_cons_unsch}, we have \begin{align*}
\label{eq2_uncon_fr_M2}
\begin{split}
  \Bignorm{\sumli |\ip{f}{g_n^*}|e^{-2\pi ic(x+\lambda_n)}T_{\lambda_n}g}_{M^1(\R)} &\Eq  \Bignorm{\sumli |\ip{f}{g_n^*}|e^{-2\pi i\lambda_nc}T_{\lambda_n}g}_{M^1(\R)}\\
  & \Le C\,\norm{f}_{M^1(\R)},
  \end{split}
\end{align*}
for some $C>0.$ Since $M^1(\R)\hookrightarrow C_0(\R)$, we obtain that $$\Bignorm{\sumli |\ip{f}{g_n^*}|\widehat{g}(\xi+c)e^{2\pi i\lambda \xi}}_{C_0(\R)} \Le C'\norm{f}_{M^1(\R)},$$
for some $C'>0.$ Consequently, we have
\begin{equation}
\label{coef_control_M1}
\sumli |\ip{f}{g_n^*}|\Le \dfrac{C'}{\widehat{g}(c)}\norm{f}_{M^1(\R)},
\end{equation}
for all $f\in M^1(\R).$
 Next, for an arbitrary $R>0$ we let $P(x)$ be an arbitrary polynomial and let $\tilde{P}(x) \in S(\R)$ be such that $\tilde{P}(x)=P(x)$ in $[-R,R]$.
     By Theorem \ref{tf_lemmas}(e), we have that $\mathcal{F}\bigparen{M^1(\R)\ast M^1(\R)}\subseteq M^1(\R)$. Consequently, we have \begin{equation}
     \label{eq2}
     \widehat{g}(x)\tilde{P}(x) \Eq \sumli \bigip{g\ast \Fc^{-1}(\tilde{P})}{g_n^*} \,e^{-2\pi i\lambda_n x }\widehat{g}(x)
    \end{equation}
     with the unconditional convergence of the series in $M^1(\R).$ Here $\mathcal{F}$ is used to denote the Fourier transform. Equation (\ref{eq2}) implies that there exists some subsequence $\set{N_j}_{j\inN}$ of $\N$ such that
     $$\widehat{g}(x)\tilde{P}(x) \Eq \lim_{j\rightarrow \infty}\sum_{n=1}^{N_j} \bigip{g\ast \Fc^{-1}(\tilde{P})}{g_n^*} \,e^{-2\pi i\lambda_n x }\widehat{g}(x) \quad $$
     pointwise almost everywhere. Combining Inequality (\ref{coef_control_M1}) with the fact that $\widehat{g}\neq 0$ almost everywhere, we see that 
    $\sumli \bigip{g\ast \Fc^{-1}(\tilde{P})}{g_n^*} e^{-2\pi i\lambda_n x }$ converges uniformly to $\tilde{P}$. However, due to the construction of $\tilde{P}$, that would imply the series $\sumli \bigip{g\ast \Fc^{-1}(\tilde{P})}{g_n^*} e^{-2\pi i\lambda_n x }$ converges uniformly to $P$ in $C([-R,R])$. Therefore, $\set{e^{-2\pi i\lambda_n x}}_{\inN}$ is complete in $C([-R,R])$ for any $R>0$.
\end{proof}

We are now ready to prove that $M^1(\R)$ does not admit any unconditional frame formed by a system of discrete translates, which consequently confirms that $M^1(\R)$ does not admit any unconditional basis formed by a system of translates by Lemma \ref{Schauder_discre_index}. We remark that the proof of the following corollary is in fact of the proof of \cite[Lemma 3.3]{OSSZ11} in disguise. Another proof of the following corollary using the perspective of Paley-Wiener space can be found in \cite[Corollary 4.33]{OU16}.

\begin{corollary}
    There does not exist any unconditional frame formed by a sequence of discrete translates for $M^1(\R)$.

    Consequently, $M^1(\R)$ does not admit an unconditional basis formed by a sequence of translates.
\end{corollary}
\begin{proof}
By Proposition \ref{effective_density_fr_M1}, it suffices to show that any discrete subset $\Lambda=\set{\lambda_n}\inN \subseteq \R$ has a finite effective density.  Let $\delta>0$ be such that $\min_{i\neq j} |\lambda_i-\lambda_j|>\delta.$ Then for any interval $(a,b]$ contained in $\R^+$, we see that $\displaystyle \max\bigparen{\frac{|(a,b] \cap  \Lambda |}{b-a}, \frac{|(a,b] \cap ( -\Lambda) |}{b-a}}< 2(\frac{b-a}{\delta})$. Therefore, we obtain that $\widetilde{D}(\Lambda)\leq \dfrac{2(b-a)}{\delta}.$ \qedhere
\end{proof}

Finally, we study the existence of unconditional frames formed by an arbitrary sequence of translates in $M^1(\R).$ The orthogonality of Rademacher system played an critical role in the proof of Theorem \ref{non_exist_discrete_unschtr}. However, this argument does not extend to $M^1(\R)$ because Rademacher functions are not in $M^1(\R).$ To obtain an ``orthogonal" sequence $\set{\psi_n}\inN$ in $M^1(\R)$, i.e. $\ip{\psi_k}{\psi_n}=\delta_{kn}$ for $k,n\in\N$, we will utilize the construction of Wilson bases. 

It was discovered by Daubechies, Jaffard and Journ\'{e} in \cite{DJJ91} that if $g\in L^2(\R)$ satisfies:
\begin{enumerate}
\setlength\itemsep{0.5em}
    \item [\textup{(a)}]$\norm{g}_{L^2(\R)}=1$, 
    \item [\textup{(b)}] $g(x)=\overline{g(-x)}$ for all $x\in \R$,
    \item [\textup{(c)}] There exists a positive constant $A$ such that $\sum_{n,k\in \Z} |\ip{f}{M_nT_{\frac{k}{2}}\psi}|^2=A\norm{f}_{L^2(\R)}$ for all $f\in L^2(\R)$,
\end{enumerate}
then 
$\set{c_nT_{\frac{k}{2}}(M_n+(-1)^{k+n}M_{-n})g}_{k\in \Z,\,n\geq 0}$, where $c_0=\frac{1}{2}$ and $c_n=\frac{\sqrt{2}}{2}$ for $n\geq 1$, is an orthonormal basis for $L^2(\R)$. Furthermore, there exists some $\phi\in C^\infty_c(\R)$ supported in $[0,1]$ that satisfies statements (a)--(c). Let $\phi_n=c_n(M_n+(-1)^nM_{-n})\phi$. Clearly, $\set{\phi_n}\inN$ satisfies the following properties 
\begin{enumerate}
\setlength\itemsep{0.5em}
    \item [\textup{(a)}]supp($\phi_n$) $\subseteq [0,1]$, 
    \item [\textup{(b)}] $\sup_n \norm{\phi_n}_{M^q(\R)}\lesssim_q \norm{\phi}_{M^1(\R)}$ for all $1\leq q \leq \infty,$
    \item [\textup{(c)}] $\ip{\phi_m}{\phi_n}=\delta_{mn}$ for all $m,n\in \N$.
    \end{enumerate}
    For the construction of such a sequence, we refer to  \cite[Corollary 8.5.5]{Gro01}.

    Now we can prove an analogue of Lemma \ref{decay_lemma_Rade} as follows. We remark that the proof follows closely to \cite[Lemma 6.3.2]{AK06}. 
    \begin{lemma} 
    \label{decay_lemma_M1}
    There exists a sequence $\set{\phi_n}\inN\subseteq C^\infty_c(\R)$ satisfying
    \begin{enumerate}
    \setlength\itemsep{0.5em}
        \item [\textup{(a)}] supp$(\phi_n)\subseteq [0,1]$ for all $n\in \N$,
\item [\textup{(b)}] $\ip{\phi_n}{\phi_m}=\delta_{mn}$,
        
        \item [\textup{(c)}] $\sup_n\norm{\phi_n}_{M^1(\R)}<\infty.$
    \end{enumerate}
     for which $\lim\limits_{n\rightarrow \infty}\ip{\phi_n}{f}=0$ for any $f\in M^p(\R)\cap L^q(\R)$ and any $(p,q)\in [1,\infty)\times [1,\infty]$. 
    \end{lemma}
    \begin{proof}
     Let $\set{\phi_n}\inN$ be the sequence that was mentioned in the paragraph right before Lemma \ref{decay_lemma_M1}. Since $\set{\phi_n}\inN$ is an orthonormal sequence in $L^2(\R)$, we see that $\lim\limits_{n\rightarrow \infty}\ip{\phi_n}{f}=0$ for any $f\in L^2(\R)$. Using that fact that $\sup_n \norm{\phi_n}_{L^\infty(\R)}\lesssim \norm{\phi}_{L^\infty(\R)}$ and the density of $L^2(\R)$ in $L^1(\R)$, we see that $\lim\limits_{n\rightarrow \infty}\ip{\phi_n}{f}=0$ for any $f\in L^1(\R)$, hence for any $f\in L^q(\R)$ for any $1\leq q\leq \infty.$ Finally, arguing similarly to Theorem \ref{decay_lemma_Rade}, we see that $\lim\limits_{n\rightarrow \infty}\ip{\phi_n }{f}=0$ for any $f\in M^p(\R).$
    \end{proof}

    The next proposition shows that if $\set{T_{\lambda_n}g}\inN$ is an unconditional frame for $M^1(\R)$ with coefficient functionals $\set{g_n^*}\inN$, then $\set{g_n^*}\inN$ must be ``rare" in the following sense.
    \begin{proposition}
    \label{no_uncon_fr_M1}
    Let $1\leq p<\infty.$
        Assume that $\schtr$ is an unconditional frame for $M^1(\R)$ with coefficient functionals $\set{g_n^*}\inN$. Then there exists some $n_0\in \N$ for which $g_{n_0}^*\notin M^p(\R)\cup L^q(\R)$ 
        for any $(p,q)\in  [1,\infty)\times [1,\infty].$
    \end{proposition}
    \begin{proof}
Suppose to the contrary that for all $n\in \N$ we have $g_n^*\in M^{q_n}(\R)$ for some $1\leq q_n<\infty$. Let $\set{\phi_n}\inN$ be the sequence described in Lemma \ref{decay_lemma_M1}. Then for any $k,N\in \N$ we have 
\begin{align*}
    \begin{split}
1&=\sumli \ip{\phi_k}{g_n^*}\ip{T_{\lambda_n}g}{\phi_k}\\
&\Le \sum_{n=1}^N \ip{\phi_k}{g_n^*}\ip{T_{\lambda_n}g}{\phi_k} + \sum_{n=N+1}^\infty \ip{\phi_k}{g_n^*}\ip{T_{\lambda_n}g}{\phi_k}\\
&\Le\sup_n\norm{\phi_n}_{M^1(\R)}\norm{g}_{M^1(\R)}\sum_{n=1}^N |\ip{\phi_k}{g_n^*}| + \sum_{n=N+1}^\infty |\ip{\phi_k}{g_n^*}||\ip{T_{\lambda_n}g}{\phi_k}|.
\end{split}
        \end{align*}
Using the estimate obtained in Equation (\ref{coef_control_M1}), we see that $|\ip{\phi_k}{g_n^*}|\leq C\sup_n\norm{\phi_k}_{M^1(\R)},$ for some $C>0$ for all $k\in \N$. On the other hand, since $M^1(\R)\hookrightarrow L^1(\R)$, we have $$\sum_{n=1}^\infty |\ip{T_{\lambda_n}g}{\phi_k}|\Le \sup_{k}\norm{\phi_k}_{L^\infty(\R)}\sum_{n=1}^\infty \norm{T_{\lambda_n}g\cdot \chi_{_{[0,1]}}}_{L^1(\R)}.$$
Choosing $N$ large enough and using Lemma \ref{decay_lemma_M1}, we obtain a contradiction. The case that $g_n^*\in L^{q_n}(\R)$ for some $1\leq q_n\leq \infty$ and all $n\in \N$ can be proved similarly. 
    \end{proof}
\begin{remark}
   We remark that the extension of Lemma \ref{decay_lemma_M1} to $M^\infty(\R)$ implies that there does not exist an unconditional frame formed by a sequence of translates for $M^1(\R).$ However, the sequence $\set{\phi_n}\inN$ that was constructed (via the construction of Wilson basis) in Lemma \ref{decay_lemma_M1} does not apply to $M^\infty (\R).$ The following example was pointed out to us by Kobayashi Masaharu. Note that $\phi_n=\sqrt{2}\cos(2\pi nx)\phi(x)$ if $n$ is even, where $\phi$ is the function mentioned in the paragraph right before Lemma \ref{decay_lemma_M1}. Let $\frac{a}{b}\in (0,1)$ be a rational number for which $|\phi(\frac{a}{b})|>0$. Let $\delta$ be the delta function. A straightforward calculation shows that $\delta\in M^\infty(\R).$ However, $|\ip{\phi_{2an}}{T_{\frac{a}{b}}\delta}|=|\sqrt{2}\phi(\frac{a}{b})|$ for all $n\in \N.$   
\end{remark}

We close this paper by showing that the result due to Odell et al.\ (\cite[Proposition 2.2]{OSSZ11}) can be extended to modulation spaces as well (see also \cite[Theorem 5.1]{DH00}).
    \begin{proposition} Let $1\leq p<\infty.$  
    \begin{enumerate}\setlength\itemsep{0.5em}
        \item [\textup{(a)}]Assume that $\set{\lambda_n}\inN$ is an uniformly discrete subset of $\R$. If $g\in M^p(\R)\cap L^1(\R)$, then $\set{g(x-\lambda_n)}\inN$ is never a Schauder frame for $M^p(\R)$ for any $1<p<\infty.$ 
     \item [\textup{(b)}] There does not exist a Schauder basis formed by a sequence of translates for $M^1(\R)$.
    \end{enumerate}
    \end{proposition}
    \begin{proof} (a)
        Suppose to the contrary that $\set{T_{\lambda_n}g}\inN$ is a Schauder frame for $M^p(\R)$ with coefficient functionals $\set{g_n^*}\inN.$ Since $\norm{T_{\lambda_n}g}_{M^p(\R)}\Eq\norm{g}_{M^p(\R)}$ for all $n\in \N$, we see that $\set{g_n^*}\inN$ converges weak$^{*}$ to $0$. Consequently, $\sup_n\norm{g_n^*}_{M^{p'}(\R)}<\infty$. Let $\set{\phi_n}\inN$ be the sequence described in Lemma \ref{decay_lemma_M1}. Then we have $|\ip{\phi_n}{g_n^*}|<\infty.$ Arguing similarly to the proof of Proposition \ref{no_uncon_fr_M1}, we obtain a contradiction.\medskip

        (b) Similarly, suppose to the contrary that $\set{T_{\lambda_n}g}\inN$ is a Schauder basis for $M^1(\R)$ with coefficient functionals $\set{g_n^*}\inN.$ Then we also have $\sup_n\norm{g_n^*}_{M^\infty(\R)}<\infty$. Using a similar argument to part (a), we obtain a contradiction. \qedhere
        
    \end{proof}
\begin{remark}
\label{final_remark}
(a) All results presented in this paper can be extended to $\R^d$ for any $d>1$ by using tensor products. For example, we outline how to extend Theorem \ref{non_exist_discrete_unschtr} (a) from $\R$ to $\R^d$ as follows. By using Lemma \ref{discrete_hilbert} step by step, we obtain $$\bignorm{\sum_{\substack{n_1,...,n_d\,\in \,\Z\\ (n_1,\dots,n_d)\neq (m_1,\dots ,m_d)}}c_{n_1,\dots,n_d} \cdot \frac{1}{m_1-n_1}\cdots\frac{1}{m_d-n_d}}_{\ell^p(\Z^d)}\Le C_p \norm{c}_{{\ell^p}(\Z^d)}$$
 for some $C_p>0$ and any $c=(c_{n_1,\dots,n_d})_{n_1,\dots,n_d\in \Z}\in \ell^p(\Z^d)$ for all $1<p<\infty.$ Next, the $\R^d$-version of Theorem \ref{painless_non_expan} states that there exists some $\phi\in C^\infty_c(\R^)$ supported in $[0,1]$ such that for any $f\in M^p(\R^d)$ there exists a sequence of scalars $(c_{kn})_{k,n\in \Z^d}\in \ell^p(\Z^{2d})$ for which 
$$f=\sum_{k,n\in \Z^d}c_{kn}M_{n}T_{\alpha k}(\underbrace{\phi \otimes\cdots \otimes \phi}_{d\text{-times}})$$
with the unconditional convergence of the series in $M^p(\R)$ for all $f\in M^p(\R^d)$ and all $1\leq p<\infty.$  Moreover, we have $\norm{f}_{M^p(\R^d)}\approx \norm{(c_{kn})_{k,n\in \Z^d}}_{\ell^p(\Z^{2d})}\approx \norm{(\ip{f}{M_nT_{\alpha k}\phi})_{k,n\in \Z^d}}_{\ell^p(\Z^{2d})}$. Here $\underbrace{\phi \otimes\cdots \otimes \phi}_{d\text{-times}}$ is the function defined on $\R^d$ by $\underbrace{(\phi \otimes\cdots \otimes \phi)}_{d\text{-times}}(x_1,\dots,x_d)=\phi(x_1)\dots\phi(x_d).$ Consequently, Theorem \ref{conv_Gabor_ext}, Lemma \ref{discrete_infinite_sum_Mpnorm}, \ref{uniform_bound_Mpnorm_Ra} and \ref{decay_lemma_Rade} can all be extended to $\R^d$ with similar arguments. It remains to extend Lemma \ref{coef_contro_Mp_unfr} to $\R^d$, then the $\R^d$-version of Theorem \ref{non_exist_discrete_unschtr} follows by a similar argument to the proof of Theorem \ref{non_exist_discrete_unschtr}. To extend Lemma \ref{coef_contro_Mp_unfr}, we consider the sequence $\set{\underbrace{R_n \otimes\cdots \otimes R_n}_{d\text{-times}}}\inN$, where $\set{R_n}\inN$ denotes the Rademacher system. Arguing similarly to the proof of $\R^1$-case, we then obtain the $\R^d$-version of Lemma \ref{coef_contro_Mp_unfr}.\medskip  

 (b) There does exist an unconditional frame formed by a sequence of translates for $M^p(\R^d)$ for any $2<p<\infty$ and any $d\geq 1.$ It was first shown by Freeman et al.\ that $L^p(\R^d)$ admits unconditional frames formed by a sequence of translates for any $p>2$ by explicitly constructing such unconditional frames. We observed that this method of construction extends to $M^p(\R^d)$ with a slight modification. Starting with the normalized unconditional basis $\set{M_nT_{k}(\chi_{_{[0,1]^d}} )}_{n,k\in \Z^d}$ for $M^p(\R^d)$, we then proceed with almost the same argument as the proof \cite[Theorem 3.2]{FOSZ14}. The only difference is that we will utilize the fact that $L^p(\R^d)\hookrightarrow M^p(\R^d)$ when $p>2$ to show the generator of the desired unconditional frame formed by a sequence of translates is in $M^p(\R^d).$\medskip

  (c) Combining results presented in this paper with Remark \ref{final_remark} (b), the only two missing pieces regarding the existence of unconditional frames and unconditional bases formed by a sequence of translates in $M^p(\R^d)$ are the answers to the following two questions. \begin{enumerate} \setlength\itemsep{0.5em}
      \item [\textup{(i)}] Assume that $\set{T_{\lambda_n}g}\inN$ is an unconditional frame for $M^1(\R)$ with coefficient functionals $\set{g_n^*}\inN$. Must $\set{g_n^*}\inN$ come solely from $ \cup_{1\leq p<\infty} M^p(\R^d)$?
      \item [\textup{(ii)}] Assume that $2<p<\infty.$ Do there exist some $g\in M^p(\R^d)$ and some $\set{\lambda_n}\inN\subseteq \R^d$ such that $\set{T_{\lambda_n}g}\inN$ forms an unconditional basis for $M^p(\R^d)$? If the answer to this question was true, then the tempered distribution $g$ that generates such a basis must be ``rare" in the sense of Theorem \ref{non_exist_discrete_unschtr} (b).
  \end{enumerate}
  
  It is known that $L^1(\R)$ does not admit an unconditional frame because $L^1(\R)$ is not topologically isomorphic to any complemented subspace of a Banach space that admits an unconditional basis (see \cite[Theorem 2.6]{CHL99} and \cite[Theorem 6.3.3]{AK06}). Here we say two Banach spaces $X$ and $Y$ are \emph{topologically isomorphic} if there exists a bijective linear operator from $X$ onto $Y$.
  However, this reasoning does not apply to $M^1(\R)$ since $M^1(\R)$ is topologically isomorphic to $\ell^1(\Z)$. 
  On the other hand, Freeman et al.\ showed in \cite[Theorem 2.1]{FOSZ14} that $L^p(\R)$ does not admit an unconditional basis formed by a sequence of translates for $p>2$ since $L^p(\R)$ is not topologically isomorphic to $\ell^p(\Z)$. Similarly, this is not the case for $M^p(\R)$ since $M^p(\R)$ is topologically isomorphic to $\ell^p(\Z^{2})$ (see \cite[Theorem 1]{FGW92}).   \qeddef
  

\end{remark}

\section{acknowledgement}
We thank Christopher Heil for his helpful comment on the manuscript. We also express our gratitude to Kobayashi Masaharu for showing us the example mentioned in Remark 3.18. Finally, we thank Nir Lev for pointing out citation errors in the Introduction.

\end{document}